\documentclass[12pt]{amsart}

\usepackage{latexsym, amsmath, amssymb, amsfonts, amscd, amsthm, mathrsfs}
\usepackage[all]{xy}

\newcommand{\Ce}{{\mathbb C}}
\renewcommand{\Re}{{\mathbb R}}
\newcommand{\Qe}{{\mathbb Q}}
\newcommand{\Ze}{{\mathbb Z}}

\newcommand{\Te}{{\mathbb T}}

\newcommand{\spn}{{\rm span}}

\newcommand{\Hom}{{\rm Hom}}
\newcommand{\End}{{\rm End}}
\newcommand{\rank}{{\rm rank}}
\newcommand{\diag}{{\rm diag}}
\newcommand{\SO}{{\rm SO}}
\newcommand{\GL}{{\rm GL}}
\newcommand{\rO}{{\rm O}}
\newcommand{\op}{{\rm op}}
\newcommand{\rL}{{\rm Lie}}
\newcommand{\tor}{{\rm tor}}

\newcommand{\cT}{{\mathcal T}}

\newcommand{\pa}{\|}

\theoremstyle{definition}
\newtheorem{theorem}{Theorem}[section]
\newtheorem{corollary}[theorem]{Corollary}
\newtheorem{lemma}[theorem]{Lemma}
\newtheorem{proposition}[theorem]{Proposition}
\newtheorem{remark}[theorem]{Remark}

\newtheorem{example}[theorem]{Example}

\begin{document}

\title[Morita equivalence of noncommutative tori]
{Strong Morita equivalence of higher-dimensional noncommutative tori. II}
\author{George A. Elliott}
\address{Department of Mathematics \\
University of Toronto \\
Toronto, Ontario, Canada~ M5S 3G3} \email{elliott@math.toronto.edu}

\author{Hanfeng Li}
\address{Department of Mathematics \\
University of Toronto \\
Toronto, Ontario, Canada~ M5S 3G3} \email{hli@fields.toronto.edu}

\date{February 7, 2005}

\subjclass[2000]{Primary 46L87; Secondary 58B34}

\begin{abstract}
We show that two $C^*$-algebraic noncommutative tori
are strongly Morita equivalent if and only if they have isomorphic
ordered $K_0$-groups and centers, extending N.~C.~Phillips's result in
the case that the algebras are simple.
This is also generalized to the twisted group $C^*$-algebras of
arbitrary finitely generated abelian groups.
\end{abstract}

\maketitle

\section{Introduction} \label{introduction:sec}

Let $n\ge 2$ and denote by $\cT_n$ the space of $n\times n$ real
skew-symmetric matrices. For each $\theta \in \cT_n$, the
corresponding $n$-dimensional ($C^*$-algebraic) noncommutative torus $A_{\theta}$ is
defined as the universal $C^*$-algebra generated by unitaries
$U_1, {\cdots}, U_n$ satisfying the relations
\begin{eqnarray*}
U_kU_j=e(\theta_{kj})U_jU_k,
\end{eqnarray*}
where $e(t)=e^{2\pi it}$. Noncommutative tori are one of the
canonical examples in noncommutative differential geometry
\cite{Rieffel90, Connes94}.

One may also consider the smooth version $A^{\infty}_{\theta}$ of
a noncommutative torus, which is the algebra of formal series
\begin{eqnarray*}
\sum c_{j_1, {\cdots}, j_n}U^{j_1}_1\cdots U^{j_n}_n
\end{eqnarray*}
where the coefficient function  $\Ze^n\ni(j_1, {\cdots}, j_n)\mapsto
c_{j_1, {\cdots}, j_n}$ belongs to the Schwartz space
$\mathcal{S}(\Ze^n)$. This is the space of smooth elements of
$A_{\theta}$ for the canonical action of $\mathbb{T}^n$ on
$A_{\theta}$.

A notion of Morita equivalence of $C^*$-algebras (as an analogue of
Morita equivalence of unital rings \cite[Chapter 6]{AF}) was
introduced by Rieffel in \cite{Rieffel74, Rieffel82}. This is now
often called strong Morita equivalence or Rieffel-Morita equivalence.
Strongly Morita equivalent $C^*$-algebras share a lot in common such
as equivalent categories of Hilbert $C^*$-modules, isomorphic
$K$-groups, etc., and hence are usually thought of as having the
same geometry.

In \cite{Schwarz98} Schwarz introduced the notion
of complete Morita equivalence of smooth noncommutative tori,
which includes strong Morita equivalence of the corresponding
$C^*$-algebras, but is stronger (as it also involves the smooth structure).
 This has
important application in M(atrix) theory \cite{Schwarz98, KS02}.

A natural question is to classify noncommutative tori and their
smooth counterparts up to the various notions of Morita
equivalence. Such results are of some interest from the point of
view of physics  \cite{CDS98, Schwarz98}. In the case $n=2$ this
was done by Rieffel \cite{Rieffel81}. In this case it does not
matter what kind of Morita equivalence one considers: there is a
(densely defined) action of the group $\GL(2, \Ze)$ on $\cT_2$,
and two matrices in $\cT_2$ give rise to Morita equivalent
noncommutative tori or smooth noncommutative tori if and only if
they are in the same orbit of this action, and also if and only if
the ordered $K_0$-groups of the algebras are isomorphic. The
higher dimensional case is much more complicated and there are
examples showing that the smooth counterparts of two strongly
Morita equivalent noncommutative tori may fail to be Morita
equivalent (as unital $\Ce$-algebras) \cite{EL} (see also
\cite{Rieffel99a, BCEN}).

After the work of Rieffel, Schwarz, and the second-named author in
\cite{Rieffel99a, Schwarz98, LiMorita} (see also \cite{TW}) it is
now known that $n$-dimensional smooth noncommutative tori are
classified up to complete Morita equivalence by the
(densely defined) $\SO(n, n|\Ze)$ action on $\cT_n$ introduced in
\cite{Rieffel99a}, as a generalization of the $\GL(2, \Ze)$ action
in the $2$-dimensional case. In \cite{EL} we showed that, in the
generic case, $n$-dimensional smooth noncommutative tori
are
classified up to Morita equivalence (as unital $\Ce$-algebras) in
the same way.

Phillips showed that two simple noncommutative tori are
strongly Morita equivalent if and only if their ordered
$K_0$-groups are isomorphic \cite[Remark 7.9]{Phillips03} (using \cite[Theorem 7.6]{Phillips03},
\cite[Theorem 1.2]{BGR},
\cite[Theorem 3.5]{Stevens}, and \cite[Lemma 3.2]{Slawny}). In this
paper we shall complete the classification of noncommutative tori
with respect to strong Morita equivalence. It turns out that in
the general (nonsimple) case one needs to know, in addition to the ordered $K_0$-group of the
algebra, only the center. (Since the center is isomorphic to $C(\Te^k)$ for some
nonnegative integer $k$, it is enough to know the dimension of $\Te^k$.
See Theorem~\ref{Morita:thm} and Remark~\ref{center:remark}.)

We shall also do more than classify the noncommutative tori (as defined above), up
to strong Morita equivalence.
 Note that $n$-dimensional noncommutative tori
are exactly the twisted group $C^*$-algebras of $\Ze^n$ (see
Subsection~\ref{group algebra:subsec}).
We shall obtain
the classification, up to strong Morita equivalence,
of the twisted group $C^*$-algebras of arbitrary finitely
generated groups.
These are exactly
the $C^*$-algebras admitting ergodic actions of compact abelian
Lie groups \cite{OPT80}.

\begin{theorem} \label{Morita:thm}
Two twisted group $C^*$-algebras of finitely generated abelian groups, in particular
two noncommutative tori,
are strongly Morita equivalent if and only
if they have isomorphic ordered $K_0$-groups and centers.
\end{theorem}

It is known that two
unital $C^*$-algebras are strongly Morita equivalent if and only if
they are Morita equivalent as unital
$\Ce$-algebras
\cite[Theorem 1.8]{Beer}.
Thus Theorem~\ref{Morita:thm} also classifies the twisted group $C^*$-algebras
of finitely generated abelian groups up to Morita equivalence as unital $\Ce$-algebras.

The opposite algebra $A^{\op}$ of a $C^*$-algebra $A$ is the
algebra $A$ with the multiplication reversed but all other
operations, including the scalar multiplication, the same. (It is
still a $C^*$-algebra.) A unital $C^*$-algebra may not be strongly
Morita equivalent to its opposite algebra (see \cite{Phillips01}
for some interesting examples), and a smooth noncommutative torus
may not be Morita equivalent to its opposite algebra (see the
discussion after Theorem 1.1 in \cite{EL}). It is a long-standing
open question whether every noncommutative torus is isomorphic to
its opposite algebra. This is now known to be the case for simple
noncommutative tori \cite[Corollary 7.8]{Phillips03}. Since the
property of admitting an ergodic action of a fixed compact group
is preserved under passing to the opposite algebra, the class of
the twisted group $C^*$-algebras of finitely generated abelian groups
is closed under the operation of taking the opposite algebra. As a
consequence of Theorem~\ref{Morita:thm}, we obtain the answer to
the easier version of the question above with strong Morita
equivalence in place of isomorphism:

\begin{corollary} \label{ME to opposite:coro}
Any twisted group $C^*$-algebra of a finitely generated abelian group is strongly Morita equivalent
to its opposite algebra.
\end{corollary}

There are two main ingredients in our proof of Theorem~\ref{Morita:thm}.
The first is that, geometrically speaking, every smooth noncommutative torus is completely Morita equivalent to
 the Cartesian product of a smooth simple noncommutative torus and
an ordinary smooth torus. This result depends on  Rieffel's
construction of finitely generated projective modules over
noncommutative tori in \cite{Rieffel88} and also on the $\SO(n,
n|\Ze)$ action on $\cT_n$ introduced by Rieffel and Schwarz in
\cite{Rieffel99a}. The second ingredient is Phillips's
classification, up to strong Morita equivalence, of simple
noncommutative tori, which as mentioned above depends on his
structure theorem for these algebras \cite[Theorem
7.6]{Phillips03}. To extend Theorem~\ref{Morita:thm} to the
twisted group $C^*$-algebras of arbitrary (countable) discrete
abelian groups, one might need to extend these results of
Phillips, Rieffel, and Schwarz to the infinite dimensional case.

This paper is organized as follows.
We review the basic definitions and facts concerning complete Morita equivalence and the twisted group algebras
in Section~\ref{preliminaries:sec}.
We prove the above-mentioned complete Morita equivalence result in Section~\ref{CME:sec} using the main result
of \cite{LiMorita}. Theorem~\ref{Morita:thm} is proved in
Section~\ref{Morita:sec} in the case of noncommutative tori.
In Section~\ref{group algebra:sec} we extend
the methods and results of Sections~\ref{CME:sec} and \ref{Morita:sec} to the twisted group algebras
of arbitrary finitely generated abelian groups, obtaining in particular
Theorem~\ref{Morita:thm} in full generality.

\section{Preliminaries} \label{preliminaries:sec}

\subsection{Complete Morita equivalence} \label{CME:subsec}

In this subsection we recall Schwarz's definition of complete Morita equivalence
and note that it passes to quotients.

We refer the reader to \cite[Sections 21 and 22]{AF} for details on algebraic Morita equivalence,
to \cite{Rieffel74, Rieffel79, Rieffel82} for strong Morita equivalence, and to
\cite{Schwarz98} and \cite[Section 7.2]{KS02}
for complete Morita equivalence.

Let $A$ and $B$ be pre-$C^*$-algebras, i.e., dense
sub-$*$-algebras of $C^*$-algebras. A strong Morita equivalence
$A$-$B$-bimodule is an $A$-$B$-bimodule ${}_AE_B$ with an
$A$-valued inner product ${}_A\left<\cdot, \cdot\right>$ and a
$B$-valued inner product $\left<\cdot, \cdot\right>_B$ satisfying
certain conditions (see \cite[Definition 6.10]{Rieffel74} for
detail; there an equivalence bimodule is called an imprimitivity
bimodule). $E$ has a norm defined by $\pa x\pa:=\pa{}_A\left<x,
x\right>\pa^{1/2}=\pa\left<x, x\right>_B\pa^{1/2}$ for $x\in E$
\cite[Proposition 3.1]{Rieffel79}. The completion of $E$ is a
strong Morita equivalence bimodule between the completions of $A$
and $B$. By \cite[Theorem 3.1, Lemma 3.1]{Rieffel79} (note that
the condition there that $A$ and $B$ be $C^*$-algebras is
unnecessary), there is a bijective correspondence between the
lattice of closed two-sided ideals $J$ of $B$ and the lattice
of closed $A$-$B$-submodules $Y$ of $E$ via
$Y=\overline{EJ}=\{y\in E: \left<x,y\right>_B\in J \mbox{ for all
} x\in E\}$ and $J=\overline{\left<E, Y\right>}$. By symmetry a similar correspondence
holds for
the lattice of closed two-sided ideals of $A$. Moveover, if
$K(J)$ is the ideal of $A$ corresponding to $Y$ (the submodule corresponding to $J$ as above),
then the
$A$-valued inner product on $E$ drops to an $A/K(J)$-valued inner
product on $E/Y$, and the $B$-valued inner product on $E$ drops to
a $B/J$-valued product on $E/Y$, so that $E/Y$ becomes a strong
Morita equivalence $A/K(Y)$-$B/J$-bimodule \cite[Corollary
3.2]{Rieffel79}. $A$ and $B$ are strongly Morita equivalent if
there exists a strong Morita equivalence $A$-$B$-bimodule. In
particular, if $A$ and $B$ are strongly Morita equivalent, then so
also are their completions.

Throughout the rest of this section we shall assume further that $A$ and $B$ are unital
and spectrally invariant, i.e., an element of $A$ (resp.~$B$) is invertible in $A$ (resp.~$B$)
if it is invertible in the completion of
$A$ (resp.~$B$). Let ${}_AE_B$ be a strong Morita equivalence $A$-$B$-bimodule.
Then ${}_A\left<E, E\right>=A$ and $\left<E, E\right>_B=B$. Furthermore,
${}_AE_B$ is an algebraic Morita equivalence $A$-$B$-bimodule (see the proof
of \cite[Theorem 1.8]{Beer}); that is,
$A=\End(E_B),\, B=\End({}_AE)$, and ${}_AE$ and $E_B$
are finitely generated projective modules and are generators
in the sense that ${}_{A}A$ and $B_B$ are direct summands
of direct sums of finitely many copies of $_{A}E$ and
$E_B$ respectively.
It is also easily checked that in the correspondence between closed two-sided
ideals $J$ of $B$ and closed $A$-$B$-submodules $Y$ of $E$ described above, one has
$Y=EJ$ (using the fact that $E_B$ is a finitely generated projective $B$-module)
and $J=\left<E, Y\right>_B$ (using $Y=EJ$ and $\left<E, E\right>_B=B$).

Suppose that a Lie algebra $L_B$ (resp.~$L_A$) acts on $B$ (resp.~$A$) as $*$-derivations.
A (Hermitian) connection on $E_B$ \cite{Connes80} is a linear map $\nabla:L_B\rightarrow \Hom_{\Ce}(E)$ satisfying the Leibniz rule:
\begin{eqnarray*}
 \nabla_X(yb)&=&(\nabla_Xy)b+y(\delta_Xb),\\
\delta_X\left<x, y\right>_B&=&\left<\nabla_Xx, y\right>_B+\left<x, \nabla_Xy\right>_B,
\end{eqnarray*}
for all $X\in L_B,\, b\in B$ and $x, y\in E$, where $\delta_X$ is the derivation of $B$ corresponding to $X$.
A connection $\nabla$ is said to have constant curvature if $[\nabla_X, \nabla_{X'}]-\nabla_{[X, X']}$ is a scalar multiplication
for all $X, X'\in L_B$. $E$ is said to be a complete Morita equivalence bimodule \cite{Schwarz98} (this is called
a gauge Morita equivalence bimodule in \cite{SW, KS02, Schwarz03}) between
$(A, L_A)$ and $(B, L_B)$ if there are constant curvature connections for ${}_AE$ and $E_B$ respectively and
a Lie algebra isomorphism from $L_B$ onto $L_A$ such that the diagram
\begin{equation*}
\xymatrix{
L_A \ar[rd] & & L_B \ar[ll] \ar[ld] \\
  & \Hom_{\Ce}(E) &
}
\end{equation*}
commutes.
Let $J,\, Y$, and $K(J)$ be as above. If a Lie subalgebra $L_{B/J}$ of $L_B$ leaves
$J$ invariant, then one checks easily that $L_{B/J}$ and $L_{A/K(J)}$ also leave $Y$ and $K(J)$
respectively invariant, where $L_{A/K(J)}$ denotes the image of $L_{B/J}$ in $L_A$, and
that the actions of $L_B$ and $L_A$ on $B$ and $A$ drop to actions of
$L_{B/J}$ and $L_{A/K(J)}$ on $B/J$ and $L_{A/K(J)}$ respectively such that
$E/Y$ is a complete Morita equivalence bimodule between
$(A/K(J), L_{A/K(L)})$ and $(B/J, L_{B/J})$.
$(A, L_A)$ and $(B, L_B)$ are said to be completely Morita equivalent if
there exists a complete Morita equivalence bimodule between them (see \cite{Schwarz99, Schwarz03}
for a more general notion called Morita equivalence of $Q$-algebras).

\subsection{Twisted group algebras} \label{group algebra:subsec}

In this subsection we shall recall basic definitions and facts about the twisted group $C^*$-algebras \cite{ZM}
and
smooth twisted group algebras of finitely generated abelian
groups.

Let $G$ be a finitely generated abelian group, and let $\sigma$ be
a $2$-cocycle on $G$ (with values in the unit circle group $\Te$),
i.e., a map $G\times G\rightarrow \Te$ such that $\sigma(g_1,
g_2)\sigma(g_1g_2, g_3)= \sigma(g_1, g_2g_3)\sigma(g_2, g_3)$ for
all $g_1, g_2, g_3\in G$. The twisted group $C^*$-algebra
$C^*(G;\sigma)$ is the universal $C^*$-algebra generated by
unitaries $u_g$ for $g\in G$ subject to the condition $u_g\cdot
u_h=\sigma(g, h)u_{gh}$. Since $G$ is finitely generated and
abelian, the dual group $\hat{G}$ is a Lie group. Because of the
universal property of $C^*(G;\sigma)$, $\hat{G}$ has a canonical
strongly continuous action $\alpha$ on $C^*(G; \sigma)$ determined
by $\alpha_x(u_g)=x(g)u_g$. This action is ergodic in the sense
that the fixed point elements are exactly the scalar multiples of
the unit. The smooth twisted group algebra $S(G; \sigma)$ is the
algebra of smooth elements of $C^*(G; \sigma)$ with respect to
this action. It is a spectrally invariant dense sub-$*$-algebra of
$C^*(G;\sigma)$. The Lie algebra $\rL(\hat{G})$ acts on $S(G;
\sigma)$ as $*$-derivations. Throughout the rest of this paper, we
shall use this Lie algebra action on $S(G; \sigma)$, and when we
talk about complete Morita equivalence between two $(S(G; \sigma),
\rL(\hat{G}))$'s we shall simply say complete Morita equivalence
between the $S(G; \sigma)$'s.

Two $2$-cocycles $\sigma$ and $\sigma'$ on $G$ are said to be
cohomologous if there is a $1$-cochain $\lambda$ on $G$, i.e., a map
$\lambda$ from $G$ to $\Te$, such that $\sigma'(g,
h)=\lambda_g\lambda_h\lambda^{-1}_{gh}\sigma(g, h)$ for all $g,
h\in G$. In this case, there is a natural $\hat{G}$-equivariant
$*$-isomorphism from $C^*(G; \sigma)$ onto $C^*(G;\sigma')$. A map
$\sigma$ from $G\times G$ to $\Te$ is said to be a bicharacter if $\sigma(g,
\cdot)$ and $\sigma(\cdot, g)$ are homomorphisms from $G$ to $\Te$
for each $g\in G$.
Clearly, every bicharacter is a $2$-cocycle. Conversely, every $2$-cocycle is cohomologous to a
bicharacter \cite[Theorem 7.1]{Kleppner}.

A bicharacter $\sigma'$ on $G$ is said to be skew-symmetric if
$\sigma(g, g)=1$ for all $g\in G$, which implies that
$\sigma(g, h)\sigma(h, g)=1$
for all $g, h\in G$.
Associated to any $2$-cocycle $\sigma$ on $G$, there is a skew-symmetric bicharacter $\sigma^*$ on $G$ defined
by $\sigma^*(g, h)=\sigma(g, h)\sigma(h, g)^{-1}$ for all $g, h\in G$.
Two $2$-cocycles $\sigma$ and $\sigma'$  on $G$ are cohomologous exactly if $\sigma^*=(\sigma')^*$ \cite[Proposition 3.2]{OPT80}.
Associated to a $2$-cocycle $\sigma$ on $G$, there is also a homomorphism $\beta_{\sigma}$ from $G$ to $\hat{G}$
defined by $\beta_{\sigma}(g)(\cdot)=\sigma^*(g, \cdot)$ for all $g\in G$.
Denote by $H_{\sigma}$ the kernel of $\beta_{\sigma}$.
The center of
$C^*(G; \sigma)$ (resp.~$S(G; \sigma)$) is the closed linear span of $u_g$, $g\in H_{\sigma}$, in $C^*(G; \sigma)$
(resp.~$S(G; \sigma)$)
and is isomorphic
to the algebra of ($\Ce$-valued) continuous (resp.~smooth) functions on $\widehat{H_{\sigma}}$ via the Fourier transform.
$\sigma$ is said to be nondegenerate if $H_{\sigma}=\{0\}$.
 It is
known that $\sigma$ is nondegenerate
exactly if $C^*(G;\sigma)$ is
simple \cite[Theorem 3.7]{Slawny}, and also exactly if $S(G;\sigma)$ is simple (see
\cite[Lemma 3.2]{BEJ84} and Theorem 13 of \cite{KS} and the remark following).

When $G=\Ze^n$ for $n\ge 0$, any $\theta\in \cT_n$ gives rise to a
skew-symmetric bicharacter
$\sigma_{\theta}$ on $\Ze^n$ defined by $\sigma_{\theta}(g, h)=e^{\pi i g\theta h^t}$.
One has $A_{\theta}=C^*(\Ze^n; \sigma_{\theta})$ and $A^{\infty}_{\theta}=S(\Ze^n; \sigma_{\theta})$
through the identification of $U_j$ with $u_{e_j}$ for $1\le j\le n$, where $U_1, {\cdots}, U_n$ are the canonical generators of
$A_{\theta}$ and $e_1, {\cdots}, e_n$ are the canonical basis elements of $\Ze^n$.
On the other hand,
any bicharacter $\sigma$ on $\Ze^n$ may be written as $\sigma(g, h)=e^{2\pi ig\Theta h^t}$ for some
$n\times n$ real matrix $\Theta$. Set $\theta=\Theta-\Theta^t$. Then $\theta$ is in $\cT_n$.
Note that $\sigma$ and $\sigma_{\theta}$ are cohomologous via the $1$-cochain $\lambda$ given by
$\lambda_g=e^{\pi ig\Theta g^t}$.
Therefore, every $2$-cocycle on $\Ze^n$ is cohomologous to $\sigma_{\theta}$ for some $\theta\in \cT_n$.
Thus, noncommutative tori and smooth noncommutative tori are exactly all the twisted group $C^*$-algebras and
smooth twisted group algebras, respectively, of torsion-free finitely generated abelian groups.

\section{Complete Morita equivalence} \label{CME:sec}

As well as the case $n\ge 2$,
we shall consider also $1$-dimensional and $0$-dimensional
noncommutative tori below, though they are actually commutative.
An element $\theta$ in $\cT_n$ for $n>0$ is said to be {\it nondegenerate} if
any element $X$ in $\Ze^n$ with $X\theta\in \Ze^n$ is $0$. As a
convention, let us say that $\theta$ is nondegenerate if $n=0$.
Then $\theta$ is nondegenerate exactly if the $2$-cocycle
$\sigma_{\theta}$ defined in Subsection~\ref{group algebra:subsec}
is nondegenerate, and also exactly if $A_{\theta}$ is simple,
and also exactly if $A^{\infty}_{\theta}$ is simple.
The main
result of this section is the following

\begin{theorem} \label{CME to tensor:thm}
For any $\theta\in \cT_n$, the smooth noncommutative torus
$A^{\infty}_{\theta}$ is completely Morita equivalent to
$A^{\infty}_{\theta'}$ for some $\theta'\in \cT_n$ such that
\begin{eqnarray*}
\theta'=\begin{pmatrix} 0 & 0 \\ 0 &  \tilde{\theta} \end{pmatrix},
\end{eqnarray*}
where $\tilde{\theta}$ belongs to $\cT_k$ for some $0\le k\le n$
and is nondegenerate.
\end{theorem}

\begin{remark} \label{CME to tensor:remark}
Note that $A^{\infty}_{\theta}$ has a natural Fr\'echet topology (see for instance
\cite[Example 3.1]{EL}). It turns out that $A^{\infty}_{\theta'}$ is the topological tensor product
\cite[Chapter 43]{Treves} of the Fr\'echet algebras $A^{\infty}_{\tilde{\theta}}$ and $C^{\infty}(\Te^{n-k})$, where
$C^{\infty}(\Te^{n-k})$ is the algebra of smooth functions on $\Te^{n-k}$.
Thus, in geometric language, Theorem~\ref{CME to tensor:thm} says that every smooth noncommutative torus
is completely Morita equivalent to the Cartesian product of a smooth simple noncommutative torus and
an ordinary smooth torus.
\end{remark}

Denote by $\rO(n, n|\Re)$ the group of linear
transformations of the space $\Re^{2n}$ preserving the quadratic
form $x_1x_{n+1}+x_2x_{n+2}+\cdots+x_nx_{2n}$, and by $\SO(n, n|\Ze)$
the subgroup of $\rO(n, n|\Re)$ consisting of matrices with integer
entries and determinant $1$.

Following \cite{Rieffel99a} let us write the elements of $\rO(n, n|\Re)$
in $2\times 2$ block form:
\begin{eqnarray*}
g=\begin{pmatrix} A & B \\ C &  D \end{pmatrix}.
\end{eqnarray*}
Then $A, B, C$, and $D$ are arbitrary $n\times n$ matrices satisfying
\begin{eqnarray*} \label{O(n,n|R):eq}
A^tC+C^tA=0=B^tD+D^tB, & A^tD+C^tB=I.
\end{eqnarray*}

The group $\SO(n, n|\Ze)$ has a partial action on $\cT_n$ \cite{Rieffel99a}, defined by
\begin{eqnarray} \label{action:eq}
g\theta=(A\theta+B)(C\theta+D)^{-1}
\end{eqnarray}
whenever $C\theta+D$ is invertible. For each $g\in \SO(n, n|\Ze)$
this action is defined on a dense open subset of $\cT_n$ (see the
discussion before \cite[Theorem]{Rieffel99a}). (The set of maximal
homeomorphisms between dense open subsets of a Hausdorff space
can be made into a group, with the natural composition law consisting of
taking ordinary composition, to the extent defined (always on a
dense open subset), and then extending to the largest open subset
on which this map is a homeomorphism onto another open subset (the
largest such open subset always exists). The map from $\SO(n, n|\Ze)$ to
the group of what might be called (maximal essential) partial
homeomorphisms of $\cT_n$ is a group homomorphism.)

Theorem~\ref{CME to tensor:thm} follows immediately from Proposition~\ref{g:prop} below and
\cite[Theorem 1.1]{LiMorita}.

\begin{proposition} \label{g:prop}
Let $\theta\in \cT_n$. Then there exists
$g\in \SO(n, n|\, \Ze)$ such that $g\theta$ is defined and
\begin{eqnarray*}
g\theta=\begin{pmatrix} 0 & 0 \\ 0 &  \tilde{\theta} \end{pmatrix},
\end{eqnarray*}
where $\tilde{\theta}$ belongs to $\cT_k$ for some $0\le k\le n$
and is nondegenerate.
\end{proposition}

\begin{example} \label{g:example}
Let $\gamma$ be a real number and
$m$ be a nonzero integer.
Consider the matrix
\begin{eqnarray*}
\theta=\begin{pmatrix} 0 & -3/m & -2/m\\ 3/m&  0 & \gamma \\ 2/m & -\gamma & 0 \end{pmatrix} \in \cT_3.
\end{eqnarray*}
Consider
$g=\begin{pmatrix} A & B \\ C &  D \end{pmatrix}\in \SO(3, 3|\, \Ze)$ where
\begin{eqnarray*}
A=\begin{pmatrix} m & 0 & 0 \\ 0 &  -2 & 3 \\  0 & -m & m\end{pmatrix}, \quad B=\begin{pmatrix} 0 & 3 & 2 \\ 0 & 0 &0 \\ 1& 0 & 0 \end{pmatrix},\\
C=\begin{pmatrix} 0 & 1 & -1 \\ 0 & 0 &0 \\ 1& 0 & 0 \end{pmatrix}, \quad D=\begin{pmatrix} 0 & 0 & 0\\ 0 &1 & 1 \\ 0 & 0 & 0  \end{pmatrix}.
\end{eqnarray*}
Then $g\theta$ is defined, and
\begin{eqnarray*}
g\theta=\begin{pmatrix} 0 & 0 &0 \\ 0 & 0& m\gamma \\ 0 & -m\gamma & 0 \end{pmatrix}.
\end{eqnarray*}
\end{example}

Let us establish some preliminary results
to prepare for the proof of Proposition~\ref{g:prop}.
Associated to the quadratic form $x_1x_{n+1}+x_2x_{n+2}+\cdots+x_nx_{2n}$ there
is a symmetric bilinear form $\left<\cdot ,\cdot \right>$ on $\Re^{2n}$ given by
\begin{eqnarray*}
\left<(x_1, {\cdots}, x_{2n}),\, (y_1, {\cdots}, y_{2n})\right>=\sum^n_{j=1}(x_jy_{n+j}+x_{n+j}y_j).
\end{eqnarray*}
The elements of $\rO(n, n|\Re)$ are exactly those linear transformations of $\Re^{2n}$ preserving
$\left<\cdot ,\cdot \right>$.
Let us say that a basis $e_1, {\cdots}, e_n, f_1, {\cdots}, f_n$ for
$\Ze^{2n}$
is compatible with the form $\left<\cdot, \cdot\right>$ if
\begin{eqnarray*}
\left<e_i, \, e_j\right>=\left<f_i, \, f_j\right>=0 \quad \mbox{ and } \quad \left<e_i, f_j\right>=\delta_{i, j}
\end{eqnarray*}
for all $1\le i, \, j\le n$. The standard basis
of $\Ze^{2n}$ is a compatible one.

\begin{lemma} \label{ON:lemma}
Let $M$ be a direct summand of $\Ze^{2n}$ such that
$M$ is isotropic, i.e.,~$\left<M, \, M\right>=0$.
Then any basis for $M$ extends to a
basis for $\Ze^{2n}$ compatible with $\left<\cdot, \cdot\right>$.
\end{lemma}
\begin{proof}
Note that the pairing between $M$ and $0^n\oplus\Ze^n$
under $\left<\cdot, \cdot\right>$ induces a homomorphism
$0^n\oplus\Ze^n\rightarrow \Hom(M, \, \Ze)$.
Denote by $W$ the kernel of this homomorphism.
Then $M+W$ is also isotropic.

Denote by $\pi$
the projection of $\Ze^{2n}$
onto $\Ze^n\oplus0^n$.
Similarly, the pairing between $\pi(M)$ and $0^n\oplus\Ze^n$
under $\left<\cdot, \cdot\right>$ induces a surjective homomorphism
$0^n\oplus\Ze^n\rightarrow \Hom(\pi(M), \, \Ze)$.
Note that the kernel of this homomorphism is also $W$.
Thus, $n=\rank(\pi(M))+\rank(W)$.

Since $\pi(M)$ is a free abelian group, we can find
a homomorphism $\psi:\pi(M)\rightarrow M$ such that
$\pi \circ \psi$ is the identity map on $\pi(M)$.
Then the restriction of $\pi$ to $\psi(\pi(M))$ is
injective, and hence $ \psi(\pi(M))\cap (0^n\oplus\Ze^n)=\psi(\pi(M))\cap \ker\pi=\{0\}$.
In particular, $\psi(\pi(M))\cap W=\{0\}$.
Consequently,
\begin{eqnarray*}
\rank(M+W)&\ge &\rank(\psi(\pi(M)))+\rank(W)\\
&=&\rank(\pi(M))+\rank(W)= n.
\end{eqnarray*}

Consider the set $\tilde{M}$ of elements of $\Ze^{2n}$
some multiple of which by a nonzero integer is in $M+W$.
Then $\tilde{M}$ is an isotropic subgroup of $\Ze^{2n}$
with  rank at least $n$.
By the elementary divisor theorem \cite[Theorem III.7.8]{Lang}, $\tilde{M}$ is a direct summand
of $\Ze^{2n}$, and any basis for $M$ extends to a basis for $\tilde{M}$.
Replacing $M$ by $\tilde{M}$, we may assume that $\rank(M)\ge n$.

The pairing between $M$ and $\Ze^{2n}$ under
$\left<\cdot, \cdot\right>$ induces a homomorphism
$\varphi: \Ze^{2n}\rightarrow \Hom(M, \, \Ze)$.
Note that for any $x\in \Ze^{2n}$, if $x/m$ is not in $\Ze^{2n}$
for every integer $m\ge 2$, then
there exists
$y\in \Ze^{2n}$ with $\left<x, \, y\right>=1$.
Using again the elementary divisor theorem, one sees
easily that $\varphi$ is surjective. Since
$M$ is contained in the kernel of $\varphi$,
we obtain
\begin{eqnarray*}
2n &= &\rank(\ker(\varphi))+\rank( \Hom(M, \, \Ze))\\
&\ge &\rank(M)+\rank( \Hom(M, \, \Ze))=2\cdot \rank(M).
\end{eqnarray*}
Therefore, $\rank(\ker(\varphi))=\rank(M)=n$. Using the elementary divisor
theorem one more time, one sees that $M$ is a direct summand of $\ker(\varphi)$.
It follows that the kernel
of $\varphi$ is exactly $M$.

Let $e_1, {\cdots}, e_n$ be a basis for $M$.
Choose $h_1, {\cdots}, h_n$ in $\Ze^{2n}$ such that
$\varphi(h_1), {\cdots}, \varphi(h_n)$ is
the dual basis of $e_1, {\cdots}, e_n$.
Then the subgroup $P$ of $\Ze^{2n}$ generated
by $h_1, {\cdots}, h_n$ maps isomorphically onto
$\Hom(M, \, \Ze)$ under $\varphi$, and hence
$\Ze^{2n}=M\oplus P$. In other words,  $e_1, {\cdots}, e_n,
h_1, {\cdots}, h_n$ is a basis for $\Ze^{2n}$.
Note that for any $x\in \Ze^{2n}$, $\left<x, \, x\right>$ is
an even integer.
Define
$f_j\in \Ze^{2n}$ inductively by
\begin{eqnarray*}
f_j=h_j-\frac{1}{2}\left<h_j, h_j\right>e_j-\sum^{j-1}_{k=1}\left<h_j, \, f_k\right>e_k.
\end{eqnarray*}
Then, clearly, $e_1, {\cdots}, e_n, f_1, {\cdots}, f_n$ is a basis for $\Ze^{2n}$ compatible with
$\left<\cdot, \cdot\right>$.
\end{proof}

\begin{remark} \label{pricipal:remark}
Let $R$ be a principal entire ring \cite[page 86]{Lang} with
characteristic not equal to $2$. Then Proposition~\ref{g:prop}
holds with $\Ze$ replaced by $R$.
\end{remark}

\begin{lemma}\label{det:lemma}
For any $n\times n$ matrix $A$ with entries in $\Ce$,
there exists a function $\zeta: \{1, {\cdots}, n\}\rightarrow \{1, -1\}$
such that
\begin{eqnarray*}
\det(A-\diag(\zeta(1), {\cdots}, \zeta(n)))\neq 0.
\end{eqnarray*}
\end{lemma}
\begin{proof}
We prove the assertion by induction on $n$.
The case $n=1$ is trivial. Suppose
that the assertion holds for $n=k$ and
$A$ is a $(k+1)\times (k+1)$ matrix.
Denote by $B$ the $k\times k$ upper left corner
of $A$. Then we can find a function
$\zeta: \{1, {\cdots}, k\}\rightarrow \{1, -1\}$
such that $\det(B-\diag(\zeta(1), {\cdots}, \zeta(k)))\neq 0$.
Define functions $\zeta^+,\, \zeta^-: \{1, {\cdots}, k+1\}\rightarrow \{1, -1\}$
extending $\zeta$ with $\zeta^{\pm}(k+1)=\pm 1$.
Observing that the matrices $A-\diag(\zeta^+(1), {\cdots}, \zeta^+(k+1))$ and
$A-\diag(\zeta^-(1), {\cdots}, \zeta^-(k+1))$ differ
at only one entry, we have
\begin{eqnarray*}
& &\det(A-\diag(\zeta^-(1), {\cdots}, \zeta^-(k+1))-\\
& &\det(A-\diag(\zeta^+(1), {\cdots}, \zeta^+(k+1)) \\
&=&2\cdot \det(B-\diag(\zeta(1), {\cdots}, \zeta(k)))\neq 0.
\end{eqnarray*}
Thus at least one of $ \zeta^+$ and $\zeta^-$ satisfies the requirement.
This finishes the induction step.
\end{proof}

\begin{lemma} \label{not inter:lemma}
Let $e_1, {\cdots}, e_n, f_1, {\cdots}, f_n$ be a
basis for $\Ze^{2n}$ compatible with $\left<\cdot, \cdot\right>$.
Let $V$ be an $n$-dimensional
isotropic linear subspace of $\Re^{2n}$, i.e., $\left<V, \, V\right>=0$.
Then we can choose $\eta_j$ from $e_j$ and $f_j$ for each $1\le j\le n$
in such a way that $\spn_{\Re}(\eta_1, {\cdots}, \eta_n)\cap V=\{0\}$.
\end{lemma}
\begin{proof}
Set $u_j=e_j+f_j$ and $v_j=e_j-f_j$ for each $j$. Set
$W_1=\spn_{\Re}(u_1, {\cdots}, u_n)$ and $W_2= \spn_{\Re}(v_1,
{\cdots}, v_n)$. With respect to the basis $u_1, {\cdots}, u_n,
v_1, {\cdots}, v_n$, the bilinear form $\left<\cdot, \cdot\right>$
on $\Re^{2n}$
gives rise to  the quadratic form
$$\Re^{2n}\ni \sum
x_ju_j+\sum y_jv_j\mapsto \sum x^2_j-\sum y^2_j.$$
Thus, the restriction of
$\left<\cdot, \cdot\right>$ to $W_2$ is negative definite.
 Since $V$ is isotropic,
$V\cap W_2=\{0\}$. Therefore, $V\subseteq \Re^{2n}=W_1\oplus W_2$
is the graph of a linear map $\varphi:W_1 \rightarrow W_2$.
Denote by $A$ the matrix of $\varphi$ with respect to the bases
$u_1, {\cdots}, u_n$ and $v_1, {\cdots}, v_n$. Choose $\zeta$ as in
Lemma~\ref{det:lemma} applied to $A$. Then $\varphi(u_1)-\zeta(1)v_1,
{\cdots}, \varphi(u_n)-\zeta(n)v_n$ are linearly independent.
Note that $u_1+\varphi(u_1), {\cdots}, u_n+\varphi(u_n)$ is a
basis for $V$. It follows easily that $\spn_{\Re}(u_1+\zeta(1)v_1,
{\cdots}, u_n+\zeta(n)v_n)\cap V=\{0\}$. Now we may just take
$\eta_j$ to be $\frac{1}{2}(u_j+\zeta(j)v_j)$ for each $j$.
\end{proof}

The next lemma is a consequence of \cite[Corollary 2.3]{TW}. For
the convenience of the reader, we give a direct proof here.

\begin{lemma} \label{SO:lemma}
Let $g=\begin{pmatrix} A & B \\ C &  D \end{pmatrix} \in \rO(n, n|\, \Ze)$.
If $C\theta+D$ is invertible for some $\theta\in \cT_n$, so that $g\theta$ is defined
in the sense of (\ref{action:eq}),
then $g$ is in $\SO(n, n|\, \Ze)$.
\end{lemma}
\begin{proof}
Suppose that $C\theta+D$ is invertible and $g$ is not in $\SO(n, n|\, \Ze)$.
Let $\varphi$ denote the linear transformation of $\Re^{2n}$
exchanging the 1st and the $(n+1)$-st coordinates. Then
$\varphi$ preserves the quadratic form $x_1x_{n+1}+x_2x_{n+2}+\cdots+x_nx_{2n}$
and has determinant $-1$. Thus the matrix $h$ corresponding
to $\varphi$ (with respect to the standard basis of $\Re^{2n}$) is
in $\rO(n, n|\, \Ze)$ but not in $\SO(n, n|\, \Ze)$.
Therefore, $hg$ is in $\SO(n, n|\, \Ze)$, and in particular
$hg$ acts on a dense open subset of $\cT_n$. Perturbing
$\theta$ slightly, we may assume then that
$(hg)(\theta)$ is defined
and that $g(\theta)$ is still defined in the sense of (\ref{action:eq}).
Set $\theta'=g(\theta)$.
Then (by matrix algebra) $\theta=g^{-1}(\theta')$ and hence (in the same way)
$hg(g^{-1}(\theta'))$ is defined in the sense of (\ref{action:eq}).
It follows (in the same way) that $h(\theta')$ is defined in the sense of (\ref{action:eq}).
But it is easy to see that $h$ does not act on any element
of $\cT_n$ in the sense of (\ref{action:eq}). Thus
we get a contradiction. Therefore, $g$ is in $\SO(n, n|\, \Ze)$.
\end{proof}

We are ready to prove Proposition~\ref{g:prop}.

\begin{proof}[Proof of Proposition~\ref{g:prop}]
Denote by $\alpha_1, {\cdots}, \alpha_n, \beta_1, {\cdots}, \beta_n$ the standard basis of $\Re^{2n}$.
Denote by $\varphi$ the linear map $0^n\oplus \Re^n\rightarrow \Re^n \oplus 0^n$
whose matrix with respect to the bases $\beta_1, {\cdots}, \beta_n$ and
$\alpha_1, {\cdots}, \alpha_n$
is $\theta$.
Since $\theta$ is skew-symmetric, the graph $V$ of $\varphi$ is an $n$-dimensional
isotropic linear subspace of $\Re^{2n}$.
Denote by $M$ the intersection of $V$ and $\Ze^{2n}$.
Using the elementary divisor theorem \cite[Theorem III.7.8]{Lang}, one sees easily that
$M$ is a direct summand of $\Ze^{2n}$. By Lemma~\ref{ON:lemma}
we can find a basis $e_1, {\cdots}, e_n, f_1, {\cdots}, f_n$
for $\Ze^{2n}$ compatible with
$\left<\cdot, \cdot\right>$ such that $f_1, {\cdots}, f_{n-k}$ is a basis
for $M$. Observing that $f_1, {\cdots}, f_{n-k}\in V$,
by Lemma~\ref{not inter:lemma} we may assume that
$\spn_{\Re}(e_1, {\cdots}, e_n)\cap V=\{0\}$.
Then $V$ is the graph of a linear map $\psi:\spn_{\Re}(f_1, {\cdots}, f_n)\rightarrow
\spn_{\Re}(e_1, {\cdots}, e_n)$. Denote by $\theta'$ the matrix of
$\psi$ with respect to the bases $f_1, {\cdots}, f_n$ and $e_1, {\cdots}, e_n$.
Since $f_1, {\cdots}, f_{n-k}\in V$, we have $\psi(f_1)={\cdots}=
\psi(f_{n-k})=0$, and hence the first $n-k$ columns of $\theta'$ are $0$.
Since $V$ is isotropic and $e_1, {\cdots}, e_n, f_1, {\cdots}, f_n$
is a basis for $\Ze^{2n}$ compatible with $\left<\cdot, \cdot\right>$, one sees easily that $\theta'$ is
skew-symmetric. Therefore, $\theta'=\begin{pmatrix} 0 & 0 \\ 0 &  \tilde{\theta} \end{pmatrix}$
for some $\tilde{\theta}\in \cT_k$.
Set $g$ to be the $2n \times 2n$ matrix such that
\begin{eqnarray*}
(e_1, {\cdots}, e_n, f_1, {\cdots}, f_n)g=(\alpha_1, {\cdots}, \alpha_n, \beta_1, {\cdots}, \beta_n).
\end{eqnarray*}
Since $g$ takes a compatible basis to
another one,
it is in $\rO(n, n|\, \Ze)$.
Using the expressions of
elements of $V$ in terms of $\theta$ and $\theta'$ respectively,
a simple calculation shows that $g\theta$ is  defined in the sense of
(\ref{action:eq})
and $g\theta=\theta'$.
By Lemma~\ref{SO:lemma}, $g$ is in $\SO(n, n|\, \Ze)$.

It remains to show that $\tilde{\theta}$ is nondegenerate.
Let $(y_{n-k+1}, {\cdots}, y_n)\in \Ze^k$ be such that
$ \tilde{\theta}(y_{n-k+1}, {\cdots}, y_n)^t$ has integral entries.
Then
\begin{eqnarray*}
\sum^n_{j=n-k+1}y_jf_j+\psi(\sum^n_{j=n-k+1}y_jf_j)\in V\cap \Ze^{2n}=M.
\end{eqnarray*}
Since $f_1, {\cdots}, f_{n-k}$ is a basis for $M$, we
get $y_{n-k+1}= {\cdots}= y_n=0$. Therefore, $\tilde{\theta}$ is nondegenerate.
\end{proof}

\section{Morita equivalence}\label{Morita:sec}

In this section we prove Theorem~\ref{Morita:thm} in the case of noncommutative tori.

We discuss first how to see whether two noncommutative tori
have isomorphic ordered $K_0$-groups.
(See \cite[Section 6]{Blackadar} for basics on ordered $K_0$-groups.)
 For any $\theta \in \cT_n$, $A_{\theta}$ has
a canonical tracial state $\tau_{\theta}$ given by the integration over the canonical action
of $\widehat{\Ze^n}$. By \cite[Lemma
2.3]{Elliott84}, all tracial states on $A_{\theta}$ induce the
same homomorphism from $K_0(A_{\theta})$ to $\Re$, which we denote by
$\omega_{\theta}$. By \cite[Theorem 3.1]{Elliott84},
$\omega_{\theta}(K_0(A_{\theta}))$ is the subgroup of $\Re$
generated by $1$ and the numbers $\sum_{\xi}(-1)^{|\xi|}\prod^{m}_{s=1}\theta_{j_{\xi(2s-1)}j_{\xi(2s)}}$ for
$1\le j_1<j_2< {\cdots}<j_{2m}\le n$, where the sum is taken over all elements $\xi$ of
the permutation group $S_{2m}$ such that $\xi(2s-1)<\xi(2s)$ for all $1\le s\le m$
and $\xi(1)<\xi(3)<\cdots<\xi(2m-1)$.
$\theta$ is said to be rational
if its entries are all rational numbers; otherwise it is said to
be nonrational. Clearly $\theta$ is rational if and only if
$\omega_{\theta}(K_0(A_{\theta}))$ has rank $1$.

\begin{proposition} \label{ordered K-0:prop}
Let $\theta_j\in \cT_{n_j}$ for $j=1,2$. Then $A_{\theta_1}$ and
$A_{\theta_2}$ have isomorphic ordered $K_0$-groups if and only if
$\omega_{\theta_2}(K_0(A_{\theta_2}))=\mu
\omega_{\theta_1}(K_0(A_{\theta_1}))$ for some real number
$\mu>0$ and either $n_1=n_2$ or $n_1+n_2=1$.
\end{proposition}
\begin{proof}
From the Pimsner-Voiculescu exact sequence \cite{PVexact} one
knows that $K_0(A_{\theta_j})$ is a free abelian group of rank $1$
or $2^{n_j-1}$ depending on whether $n_j=0$ or $n_j>0$.

We prove first the ``only if'' part. Comparing the ranks of the
$K_0$-groups we see that either $n_1=n_2$ or $n_1+n_2=1$ (i.e., one of
$n_1$ and $n_2$ is $0$ and the other is $1$). Every unital
$C^*$-algebra admitting an ergodic action of a compact abelian
group is nuclear \cite[Lemma 6.2]{OPT80}, \cite[Proposition
3.1]{DLRZ02}. Since $A_{\theta_j}$ admits an ergodic action of
$\Te^{n_j}$, it is nuclear and hence is exact. For a unital
$C^*$-algebra $A$, a state on the scaled ordered $K_0$-group
$(K_0(A)_+, K_0(A), [1_A])$ is a positive unital homomorphism from
$(K_0(A)_+, K_0(A), [1_A])$ to $(\Re_+, \Re, 1)$. When $A$ is
exact, every state on $(K_0(A)_+, K_0(A), [1_A])$ comes from a
tracial state on $A$ \cite[Corollary 3.4]{BR}, \cite[Theorem 9.2]{HT}. Therefore,
$(K_0(A_{\theta_j})_+, K_0(A_{\theta_j}), [1_{A_{\theta_j}}])$ has
a unique state, which is exactly $\omega_{\theta_j}$. Take an order
isomorphism $\psi$ from $(K_0(A_{\theta_1})_+, K_0(A_{\theta_1}))$
onto $(K_0(A_{\theta_2})_+, K_0(A_{\theta_2}))$. Then
$\omega_{\theta_2}\circ \psi$ is a nontrivial positive homomorphism from
$(K_0(A_{\theta_1})_+, K_0(A_{\theta_1}))$ to $(\Re_+, \Re)$. Set
$\mu$ equal to the value of $[1_{A_{\theta_1}}]$ under
$\omega_{\theta_2}\circ \psi$. Then $\mu>0$ and
$\frac{1}{\mu}(\omega_{\theta_2}\circ \psi)$ is a state on
$(K_0(A_{\theta_1})_+, K_0(A_{\theta_1}), [1_{A_{\theta_1}}])$.
Consequently, $\frac{1}{\mu}(\omega_{\theta_2}\circ \psi)=\omega_{\theta_1}$.
Evaluating both sides on $K_0(A_{\theta_1})$ we get
$\omega_{\theta_2}(K_0(A_{\theta_2}))=\mu \omega_{\theta_1}(K_0(A_{\theta_1}))$.

Next we prove the ``if'' part. Note that
$\omega_{\theta_j}(K_0(A_{\theta_j}))$ is a torsion-free finitely
generated abelian group, and hence is a free abelian group. Taking
a lifting of $\omega_{\theta_j}(K_0(A_{\theta_j}))$ in
$K_0(A_{\theta_j})$ and identifying this lifting with
$\omega_{\theta_j}(K_0(A_{\theta_j}))$, we may assume that
$K_0(A_{\theta_j})=\ker(\omega_{\theta_j})\oplus
\omega_{\theta_j}(K_0(A_{\theta_j}))$ and that $\omega_{\theta_j}$
is exactly the projection onto the second summand. Now we need to
distinguish the cases $\theta_j$ is rational or nonrational.
Suppose that both $\theta_1$ and $\theta_2$ are nonrational. Then
$n_1=n_2$. By \cite[Theorem 6.1]{Rieffel88},
$(K_0(A_{\theta_j}))_+$ consists of exactly the elements of
$K_0(A_{\theta_j})$ on which $\omega_{\theta_j}$ is strictly
positive, together with $0$. The multiplication by $\mu$ is an
order isomorphism from $\omega_{\theta_1}(K_0(A_{\theta_1}))$ onto
$\omega_{\theta_2}(K_0(A_{\theta_2}))$.
 Then
$\ker(\omega_{\theta_1})$ and $\ker(\omega_{\theta_2})$ have the
same rank. Taking any isomorphism from $\ker(\omega_{\theta_1})$
onto $\ker(\omega_{\theta_2})$, we get an order isomorphism from
$(K_0(A_{\theta_1})_+, K_0(A_{\theta_1}))$ onto
$(K_0(A_{\theta_2})_+, K_0(A_{\theta_2}))$, as desired.

Now assume
that at least one of $\theta_1$ and $\theta_2$ is rational.
Comparing the ranks of $\omega_{\theta_1}(K_0(A_{\theta_1}))$ and
$\omega_{\theta_2}(K_0(A_{\theta_2}))$ we see that both $\theta_1$
and $\theta_2$ are rational. Note that the partial action of
$\SO(n, n|\Ze)$ on $\cT_n$ preserves rationality. Thus, if we
take $\theta=\theta_j$ and $n=n_j$ in Proposition~\ref{g:prop},
then $\tilde{\theta}$
given in Proposition~\ref{g:prop} must be rational and hence
$k$ given there must be $0$. By \cite[Theorem 1.1]{LiMorita},
$A^{\infty}_{\theta_j}$ is completely Morita equivalent to
$A^{\infty}_{0_{n\times n}}$, where $0_{n\times n}$ is the zero
$n\times n$ matrix. Consequently, $A_{\theta_1}$ and
$A_{\theta_2}$ are Morita equivalent and hence have isomorphic ordered $K_0$-groups.
This finishes the proof of
Proposition~\ref{ordered K-0:prop}.
\end{proof}

The ``only if'' part of Proposition~\ref{ordered K-0:prop}
and the proof of the ``if'' part
show

\begin{corollary} \label{rational:coro}
Let $\theta\in \cT_n$. Then the following statements are equivalent:
\begin{enumerate}
\item $\theta$ is rational; \item $A^{\infty}_{\theta}$ is
completely Morita equivalent to $A^{\infty}_{0_{n\times n}}$,
where $0_{n\times n}$ is the zero $n\times n$ matrix; \item
$A_{\theta}$ is Morita equivalent to $A_{0_{n\times n}}=C(\Te^n)$, where $C(\Te^n)$ is the algebra of continuous functions on $\Te^n$;
\item $\omega_{\theta}(K_0(A_{\theta}))$ has rank $1$.
\end{enumerate}
\end{corollary}

\begin{lemma} \label{K-0 for simple:lemma}
Let $\theta_1$ and $\theta_2$ in $\cT_n$ be such that
\begin{eqnarray*}
\theta_1=\begin{pmatrix} 0 & 0 \\ 0 &  \tilde{\theta}_1 \end{pmatrix}\quad  \mbox{ and } \quad
\theta_2=\begin{pmatrix} 0 & 0 \\ 0 &  \tilde{\theta}_2 \end{pmatrix},
\end{eqnarray*}
with $\tilde{\theta}_1, \tilde{\theta}_2\in \cT_k$. If the
ordered $K_0$-groups of $A_{\theta_1}$ and $A_{\theta_2}$ are
isomorphic, then so also are those of
$A_{\tilde{\theta}_1}$ and $A_{\tilde{\theta}_2}$.
\end{lemma}
\begin{proof}
Note that $A_{\theta_j}=A_{\tilde{\theta}_j}\otimes C(\Te^{n-k})$.
Taking the evaluation at any point of $\Te^{n-k}$ we get a unital
$*$-homomorphism $\varphi_j$ from $A_{\theta_j}$ to
$A_{\tilde{\theta}_j}$. Denote by $(\varphi_j)_*$ the induced
homomorphism from $K_0(A_{\theta_j})$ to
$K_0(A_{\tilde{\theta}_j})$. Then $\omega_{\tilde{\theta}_j}\circ
(\varphi_j)_*$ is exactly $\omega_{\theta_j}$. Using the embedding
$A_{\tilde{\theta}_j}\hookrightarrow A_{\tilde{\theta}_j}\otimes
C(\Te^{n-k})= A_{\theta_j}$ one sees that
$(\varphi_j)_*$ is surjective.
Thus,
$$\omega_{\theta_j}(K_0(A_{\theta_j}))=(\omega_{\tilde{\theta}_j}\circ (\varphi_j)_*)(K_0(A_{\theta_j}))
=\omega_{\tilde{\theta}_j}(K_0(A_{\tilde{\theta}_j})).$$ Now
Lemma~\ref{K-0 for simple:lemma} follows from
Proposition~\ref{ordered K-0:prop}.
\end{proof}

As we mentioned in Subsection~\ref{group algebra:subsec}, the center of $A_{\theta}$ is
isomorphic to the algebra of continuous functions on $\widehat{H_{\sigma_{\theta}}}$, and hence
depends only on the rank of $H_{\sigma_{\theta}}$, which can be calculated from $\theta$ arithmetically.

We are ready to prove Theorem~\ref{Morita:thm} in the case of noncommutative tori.

\begin{proof}[Proof of Theorem~\ref{Morita:thm} in the case of noncommutative tori]
The ``only if'' part follows from the fact that Morita equivalence between
unital algebras (or rings) preserves both the
ordered $K_0$-group and the center \cite[Proposition 21.10]{AF}.
Consider the ``if'' part.
Suppose that $A_{\theta_1}$ and $A_{\theta_2}$ have isomorphic ordered
$K_0$-groups and centers. By the ``only if'' part and Theorem~\ref{CME to tensor:thm},
we may assume that
\begin{eqnarray*}
\theta_j=\begin{pmatrix} 0 & 0 \\ 0 &  \tilde{\theta}_j \end{pmatrix}
\end{eqnarray*}
for some nondegenerate $\tilde{\theta}_j \in \cT_{k_j}$.
Say that $\theta_j$ is in $\cT_{n_j}$.
Then $A_{\theta_j}\cong C(\Te^{n_j-k_j})\otimes A_{\tilde{\theta}_j}$.
From the Pimsner-Voiculescu exact sequence \cite{PVexact} one
knows (as recalled above) that $K_0(A_{\theta_j})$ is a free abelian group of rank $1$
or $2^{n_j-1}$, depending on whether $n_j=0$ or $n_j>0$. It follows that
either $n_1=n_2$, or
$n_1+n_2=1$. Note that $A_{\theta}=\Ce$ if $n_j=0$ and
$A_{\theta}=C(\Te)$ if $n_j=1$. Therefore we must have $n_1=n_2$. Also
note that the center of $A_{\theta_j}$ is isomorphic to
$C(\Te^{n_j-k_j})$. Thus $k_1=k_2$. By Lemma~\ref{K-0 for
simple:lemma}, $A_{\tilde{\theta}_1}$ and $A_{\tilde{\theta}_2}$
have isomorphic ordered $K_0$-groups. Since $\tilde{\theta}_1$ and
$\tilde{\theta}_2$ are nondegenerate, both
$A_{\tilde{\theta}_1}$ and $A_{\tilde{\theta}_2}$ are simple.
Phillips has shown that simple noncommutative tori are classified
up to strong Morita equivalence by their ordered $K_0$-groups
\cite[Remark 7.9]{Phillips03}. Therefore, $A_{\tilde{\theta}_1}$
and $A_{\tilde{\theta}_2}$ are strongly Morita equivalent.
Consequently, $A_{\theta_1}$ and $A_{\theta_2}$ are strongly
Morita equivalent. This finishes the proof of
Theorem~\ref{Morita:thm} in the case of noncommutative tori.
\end{proof}

\begin{remark} \label{center:remark}
In the case of simple noncommutative tori, note that the result of Phillips in \cite{Phillips03}
refers only to the ordered $K_0$-group, not the center.
Theorem~\ref{Morita:thm} in fact also does this as the center is the scalars in this case.
If we consider only
noncommutative tori of dimension $2$ or $3$, then
Theorem~\ref{Morita:thm} holds without mentioning the
centers, since in these cases the dimension of the center is determined
by the ordered $K_0$-group.
The reason in the case of $2$-dimensional noncommutative
tori is that in this case the dimension of the center of
$A_{\theta}$ is either $2$ or $0$, depending as
$\omega_{\theta}(K_0(A_{\theta}))$ has rank $1$ or $2$, as is easily
seen from the arithmetical description of $\omega_{\theta}(K_0(A_{\theta}))$
given above.
The reason in the case of $3$-dimensional
noncommutative tori is that in this case the center of
$A_{\theta}$ has dimension $3$, $1$, or $0$, depending as
$\omega_{\theta}(K_0(A_{\theta}))$ has rank $1$, $2$, or at least
$3$. However, Example~\ref{need center:remark} below shows
that if we consider $n$-dimensional noncommutative tori for a
fixed $n\ge 4$, then Theorem~\ref{Morita:thm} does not hold any
longer without keeping track of the centers.
\end{remark}

\begin{example} \label{need center:remark}
Let $\gamma$ be a real algebraic integer of degree $2$ (for
example, $\sqrt{2}$). Set
\begin{eqnarray*}
\theta_1=\begin{pmatrix} 0 & \gamma & 0 & 0\\ -\gamma & 0 & 0 & 0\\
0 & 0 & 0 & 0\\ 0 & 0 & 0 & 0\end{pmatrix} \quad \mbox{ and }
\quad \theta_2=\begin{pmatrix} 0 & \gamma & 0 & 0\\ -\gamma & 0 &
 0 & 0\\0 & 0 & 0&  \gamma\\ 0 & 0&  -\gamma & 0\end{pmatrix}.
\end{eqnarray*}
Then by the arithmetical description of the range of the trace on $K_0$ given above,
$\omega_{\theta_1}(K_0(A_{\theta_1}))=\omega_{\theta_2}(K_0(A_{\theta_2}))$,
and hence $A_{\theta_1}$ and $A_{\theta_2}$ have isomorphic
ordered $K_0$-groups by Proposition~\ref{ordered K-0:prop}. But
the center of $A_{\theta_1}$ has dimension $2$, while that of
$A_{\theta_2}$ has dimension $0$.
\end{example}

\section{Twisted group algebras of finitely generated abelian groups}
\label{group algebra:sec}

In this section we extend the results of Sections~\ref{CME:sec}
and \ref{Morita:sec} to the twisted group algebras of arbitrary finitely
generated abelian groups.

Denote by $G_{\tor}$ the torsion subgroup of a finitely generated abelian group $G$.
The rank of $G$ is the dimension of the $\Qe$-vector space $G\otimes_{\Ze}\Qe$.
Recall that a $2$-cocycle on $G$ is nondegenerate if the subgroup $H_{\sigma}$ defined in Subsection~\ref{preliminaries:sec}
is $\{0\}$.
The following result is a generalization of Theorem~\ref{CME to tensor:thm}.

\begin{theorem} \label{CME to tensor2:thm}
Let $\sigma$ be a $2$-cocycle on a finitely generated abelian group $G$. Then
there exist a finitely generated abelian group $G'=G'_1\oplus G'_2$ and
a skew-symmetric bicharacter $\sigma'$ on $G'$ such that
$\rank(G)=\rank(G')$,
$G'_2$ is torsion-free,
$\sigma'(G'_1, G')=1$, the restriction of $\sigma'$ to $G'_2$ is nondegenerate,
and $S(G'; \sigma')$ is completely Morita equivalent to $S(G; \sigma)$.
\end{theorem}
\begin{proof}
Since $G$ is finitely generated and abelian, we can find a nonnegative integer $n$ and
a surjective homomorphism $\psi$ from $\Ze^n$ to $G$. Then the pull-back $\psi^*(\sigma)$ of $\sigma$ under
$\psi$ is a $2$-cocycle on $\Ze^n$. In Subsection~\ref{group algebra:subsec}
we noticed that there exists an element $\theta$ of $\cT_n$ such that
$\psi^*(\sigma)$ is cohomologous to $\sigma_{\theta}$ via a $1$-cochain $\lambda$ on $\Ze^n$.
By Theorem~\ref{CME to tensor:thm}, $A^{\infty}_{\theta}$ is
completely Morita equivalent to $A^{\infty}_{\theta'}$ for some $\theta' \in \cT_n$ of the form
\begin{eqnarray*}
\theta'=\begin{pmatrix} 0 & 0 \\ 0 &  \tilde{\theta} \end{pmatrix},
\end{eqnarray*}
where $\tilde{\theta}\in \cT_k$ for some $0\le k\le n$
and is nondegenerate. Let $_{A^{\infty}_{\theta'}}E_{A^{\infty}_{\theta}}$ be a
complete Morita equivalence $A^{\infty}_{\theta'}$-$A^{\infty}_{\theta}$-bimodule
with constant-curvature connections on $E_{A^{\infty}_{\theta}}$ and $_{A^{\infty}_{\theta'}}E$ and
a Lie algebra isomorphism $\phi:\rL(\widehat{\Ze^n})\rightarrow \rL(\widehat{\Ze^n})$ as in Subsection~\ref{CME:subsec}.

In Subsection~\ref{CME:subsec} we mentioned that complete Morita equivalence passes to
quotient algebras (with the actions of certain Lie subalgebras). We shall find
an ideal $J$ of $A^{\infty}_{\theta}$ such that $A^{\infty}_{\theta}/J\cong S(G; \sigma)$ and
identify the corresponding ideal $K(J)$ of $A^{\infty}_{\theta'}$. Then
$A^{\infty}_{\theta'}/K(J)$ is completely Morita equivalent to $S(G; \sigma)$.

Note that there is a $*$-homomorphism $\psi_*$ from $A_{\theta}=C^*(\Ze^n; \sigma_{\theta})$ onto
$C^*(G; \sigma)$ sending $u_g$ to $\lambda_gu_{\psi(g)}$ for all $g\in \Ze^n$.
From the universal property of $C^*(G; \sigma)$ one sees that $\ker(\psi_*)$ is the closed
ideal of $C^*(\Ze^n; \sigma_{\theta})$ generated by $u_{g_j}-\lambda_{g_j}\lambda^{-1}_{0}, 1\le j\le m$ for any basis $g_1, {\cdots}, g_m$
for $\ker(\psi)$.  Denote $u_{g_j}-\lambda_{g_j}\lambda^{-1}_{0}$ by $v_j$.
Observing that $v_j$ is in the center of $C^*(G; \sigma_{\theta})$, we see that
$\ker(\psi_*)$ is the closure of $\sum^m_{j=1}
v_jC^*(\Ze^n; \sigma_{\theta})$ in
$C^*(\Ze^n; \sigma_{\theta})$. Denote by $J$ the intersection of $S(\Ze^n; \sigma_{\theta})$
and $\ker(\psi_*)$. Then $J$ is a closed two-sided ideal of $S(\Ze^n; \sigma_{\theta})$
and is the closure of $\sum^m_{j=1}
v_jS(\Ze^n; \sigma_{\theta})$ in
$S(\Ze^n; \sigma_{\theta})$.
Using $\psi$ we may identify
$\hat{G}$ with a closed subgroup of $\widehat{\Ze^n}$,
and hence identify $\rL(\hat{G})$ with a
Lie subalgebra of $\rL(\widehat{\Ze^n})$.
Clearly, $\psi_*$ sends $A^{\infty}_{\theta}=S(\Ze^n;\sigma_{\theta})$ into $S(G; \sigma)$.
It follows that $\rL(\hat{G})$ preserves $J$.
It is easily checked that elements of $S(G;\sigma)$ are of the form $\sum_{g\in G}c_gu_g$
for the coefficient function $G\ni g\mapsto c_g$ belonging to the Schwarz space $\mathcal{S}(G)$
(cf. the proof of the corollary on page 468 of \cite{MW}).
Thus $\psi_*$ sends $S(\Ze^n;\sigma_{\theta})$ onto $S(G; \sigma)$.
Consequently, $A^{\infty}_{\theta}/J$ and
$S(G; \sigma)$ are isomorphic as pre-$C^*$-algebras, in a way compatible
with the actions of $\rL(\hat{G})$.

Now we need to find the corresponding ideal $K(J)$ of $A^{\infty}_{\theta'}$.
There is a $*$-isomorphism $\varphi$ from the center of
$A^{\infty}_{\theta}$ onto that of $A^{\infty}_{\theta'}$
determined by $\varphi(a)x=xa$ for all $a\in A^{\infty}_{\theta}$
and $x\in E$.
Clearly, the corresponding closed $A^{\infty}_{\theta'}$-$A^{\infty}_{\theta}$-submodule $Y$ of
$E$ is the closure of $\sum^m_{j=1}Ev_j$ in $E$, and
$K(J)$ is the closed two-sided ideal of $A^{\infty}_{\theta'}$ generated by $\varphi(v_j), 1\le j\le m$.
Recall that the center of $A^{\infty}_{\theta}=S(G; \sigma_{\theta})$ is the closed linear span
of $u_g$, $g\in H_{\sigma_{\theta}}$, in $S(G; \sigma_{\theta})$ (see Subsection~\ref{group algebra:subsec}).
Moreover, the eigenvectors of the restriction of
the action
of $\rL(\widehat{\Ze^n})$ on $S(\Ze^n;\sigma_{\theta})$ to
the center are exactly the scalar multiples of the $u_g$'s for $g\in
H_{\sigma_{\theta}}$.
Therefore, for any $g\in H_{\sigma_{\theta}}$, up to a scalar
multiple, $\varphi(u_g)$ is equal to $u_{g'}$ for some
$g'\in H_{\sigma_{\theta'}}$. The map $\eta$ from
$H_{\sigma_{\theta}}$ to $H_{\sigma_{\theta'}}$ sending $g$ to
$g'$ is easily seen to be an isomorphism. Consequently, $K(J)$ is the
closed two-sided ideal of $A^{\infty}_{\theta'}$ generated by
$u_{\varphi(g_j)}-\gamma_j, 1\le j\le m$, where $\gamma_j$ is
a certain element of $\Te$. Note that $\ker(\psi)\subseteq
H_{\sigma_{\theta}}$ and $\Ze^{n-k}\oplus
0^k=H_{\sigma_{\theta'}}$. Set $G'=\Ze^n/\eta(\ker(\psi))$. Then
$G'=G'_1\oplus G'_2$ for
$G'_1=H_{\sigma_{\theta'}}/\eta(\ker(\psi))$ and $G'_2=\Ze^k$.
Set $\sigma'$ equal to the skew-symmetric
bicharacter on $G'$ such that $\sigma'(G'_1, G')=1$ and the
restrictions of $\sigma'$ and $\sigma_{\theta'}$ to $G'_2$
coincide. Then the pull-back of $\sigma'$ under the quotient map
$\psi'$ from $\Ze^n$ to $G'$ is exactly $\sigma_{\theta'}$, and
the  restriction of $\sigma'$ to $\Ze^k$ is nondegenerate.
Note
that for any basis $h_1, {\cdots}, h_{n-k}$ for $ \Ze^{n-k}$ and any
$\mu_1, {\cdots}, \mu_{n-k}\in \Te$ there is a $*$-homomorphism from
$C^*(\Ze^n; \sigma_{\theta'})=A_{\theta'}$ onto $C^*(G'; \sigma')$
sending $u_{h_j}$ to $\mu_ju_{\psi(h_j)}$ for all $1\le j\le k$
and sending $u_h$ to $u_{\psi(h)}$ for all $h\in \Ze^k$. Using the
elementary divisor theorem \cite[Theorem III.7.8]{Lang}, we can
choose suitable $\mu_1, {\cdots}, \mu_{n-k}$ such that the kernel
of the above $*$-homomorphism
is exactly the closed two-sided ideal of $C^*(\Ze^n; \sigma_{\theta'})$
generated by $u_{\varphi(g_j)}-\gamma_j, 1\le j\le m$. Then
$K(J)$ is the intersection of this ideal and
$S(\Ze^n;\sigma_{\theta'})=A^{\infty}_{\theta'}$. As in the last
paragraph, we may identify $\rL(\widehat{G'})$ with a Lie
subalgebra of $\rL(\Ze^n)$, and $K(J)$ is invariant under the
action of $\rL(\widehat{G'})$. Furthermore,
$A^{\infty}_{\theta'}/K(J)$ and $S(G';\sigma')$ are isomorphic as
pre-$C^*$-algebras, in a way compatible with the actions of
$\rL(\widehat{G'})$. Observing that $\rL(\hat{G})$
(resp.~$\rL(\widehat{G'})$) consists of exactly those elements of
$\rL(\widehat{\Ze^n})$ (resp.~$\rL(\widehat{\Ze^n})$) acting trivially on the center
of $A^{\infty}_{\theta}$ (resp.~$A^{\infty}_{\theta'}$), we see
that $\rL(\hat{G})$ is sent onto $\rL(\widehat{G'})$ under the Lie
algebra isomorphism $\phi:\rL(\widehat{\Ze^n})\rightarrow
\rL(\widehat{\Ze^n})$. Therefore, $S(G';\sigma')$ and $S(G;\sigma)$
are completely Morita equivalent.
From $\rank(G)+\rank(\ker(\psi))=\rank(\Ze^n)$ and
$\rank(G')+\rank(\eta(\ker(\psi)))=\rank(\Ze^n)$ we obtain
$\rank(G)=\rank(G')$.
This finishes the proof of
Theorem~\ref{CME to tensor2:thm}.
\end{proof}

Recall the group $H_{\sigma}$ and the skew-symmetric bicharacter $\sigma^*$
defined in Subsection~\ref{group algebra:subsec}
for a $2$-cocycle $\sigma$ on a finitely generated abelian group $G$.

\begin{remark} \label{CME to tensor2:remark}
We indicate briefly another proof of Theorem~\ref{CME to
tensor2:thm}. Note that $\sigma^*$ induces a skew-symmetric
bicharacter $\tau^*$ on $G'':=G/(H_{\sigma})_{\tor}$. Then
$\tau^*=(\sigma'')^*$ for some $2$-cocycle $\sigma''$ on $G''$.
Since $\widehat{G''}$ is a subgroup of $\hat{G}$, it also acts on
$C^*(G;\sigma)$. One checks easily that $C^*(G;\sigma)$ and the
direct sum of $|(H_{\sigma})_{\tor}|$ many copies of
$C^*(G'';\sigma'')$ are isomorphic in a way compatible with the
actions of $\widehat{G''}$. Using the proof of
\cite[Proposition]{Rosenberg}, by induction on $|G''_{\tor}|$, one
can show that there is a free (abelian) subgroup $G'$ of $G''$ of the same
rank as $G$ such that $C^*(G'';\sigma'')$ is isomorphic to
$M_m(C^*(G';\sigma'))$ for some $m$ and $S(G'';\sigma'')$ is
completely Morita equivalent to $S(G';\sigma')$, where $\sigma'$
is the restriction of $\sigma''$ to $G'$. Consequently,
$C^*(G;\sigma)$ is isomorphic to the direct sum of
$|(H_{\sigma})_{\tor}|$ many copies of
$M_m(C^*(G';\sigma'))$ and
$S(G;\sigma)$ is completely Morita equivalent to the direct sum of
$|(H_{\sigma})_{\tor}|$ many copies of $S(G';\sigma')$. Then
Theorem~\ref{CME to tensor2:thm} follows from Theorem~\ref{CME to
tensor:thm}. However, the earlier detailed proof does not use
induction and is more likely to be generalizable to the case of
arbitrary (countable) discrete abelian groups.
\end{remark}

\begin{lemma} \label{size:lemma}
In Theorem~\ref{CME to tensor2:thm}
one has $|(H_{\sigma})_{\tor}|=|(G'_1)_{\tor}|$ and
$\rank(H_{\sigma})=\rank(G'_1)$.
\end{lemma}
\begin{proof}
Recall that the center of $C^*(G; \sigma)$ is isomorphic to $C(\widehat{H_{\sigma}})$.
Note that $|(H_{\sigma})_{\tor}|$ and $\rank(H_{\sigma})$ are
the number of connected components and the dimension of
$\widehat{H_{\sigma}}$ respectively.
Since Morita equivalence
preserves centers \cite[Proposition 21.10]{AF}, we obtain $|(H_{\sigma})_{\tor}|=|(H_{\sigma'})_{\tor}|=|(G'_1)_{\tor}|$
and $\rank(H_{\sigma})=\rank(H_{\sigma'})=\rank(G'_1)$.
\end{proof}

From the Pimsner-Voiculescu exact sequence \cite{PVexact} one
knows that the $K_0$-group (resp.~$K_1$-group) of an $n$-dimensional noncommutative
torus is a free abelian group of rank
$2^{n-1}$ or $1$ (resp.~$2^{n-1}$ or $0$) depending as $n>0$ or $n=0$.
Since strong Morita equivalence between $C^*$-algebras
preserves $K$-groups \cite[Theorem 1.2]{BGR}, \cite{Exel}, we get

\begin{corollary} \label{K-groups:coro}
Let $\sigma$ be a $2$-cocycle on a finitely generated abelian group $G$. Then
the $K_0$-group (resp.~$K_1$-group) of $C^*(G;\sigma)$ is a free abelian group with rank
$|(H_{\sigma})_{\tor}|\cdot 2^{\rank(G)-1}$ or $|(H_{\sigma})_{\tor}|$
(resp.~$|(H_{\sigma})_{\tor}|\cdot 2^{\rank(G)-1}$ or $0$)
depending as $\rank(G)>0$ or $\rank(G)=0$.
\end{corollary}

For any unital $C^*$-algebra $A$, denote by $T(A)_{K_0}$ the set of all
homomorphisms from $K_0(A)$ to $\Re$ induced by tracial states of
$A$. Then $T(A)_{K_0}$ equipped with the topology of pointwise convergence
is a compact convex set in a Hausdorff locally convex topological vector space.

\begin{lemma} \label{states:lemma}
Let $\sigma$ be a $2$-cocycle on a finitely generated abelian group $G$.
Then $T(C^*(G; \sigma))_{K_0}$ is a simplex of dimension $|(H_{\sigma})_{\tor}|-1$.
The images of $K_0(C^*(G;\sigma))$ are the same under all the vertices of $T(C^*(G; \sigma))_{K_0}$.
\end{lemma}
\begin{proof} 
Recall the skew-symmetric bicharacter $\sigma^*$ on $G$ defined in
Subsection~\ref{group algebra:subsec}. It induces a skew-symmetric bicharacter $\tau^*$
on $G/H_{\sigma}$. Then $\tau^*=(\sigma')^*$ for some $2$-cocycle $\sigma'$ on $G/H_{\sigma}$.
Note that $\sigma'$ is nondegenerate.
Thus $C^*(G/H_{\sigma}; \sigma')$ has a unique tracial state $\varphi$ \cite[Lemma 3.2]{Slawny}.
For any closed ideal $I_t$ of $C^*(G; \sigma)$ generated by a maximal ideal $t$ of
the center of $ C^*(G; \sigma)$, one has $C^*(G; \sigma)/I_t\cong C^*(G/H_{\sigma}; \sigma')$.

Any extremal point of $T(C^*(G; \sigma))_{K_0}$ is induced by an extremal tracial state of
$C^*(G; \sigma)$. An argument similar to that in the proof of \cite[Lemma 2.2]{Elliott84}
shows that every extremal tracial state of
$C^*(G; \sigma)$ factors through $C^*(G; \sigma)/I_t$ for some maximal ideal $t$ of
the center (and hence must be the pull-back $\varphi_t$ of the unique tracial state of $C^*(G; \sigma)/I_t$),
and also that the map from $\widehat{H_{\sigma}}$ to $T(C^*(G; \sigma))_{K_0}$ sending
$t$ to the homomorphism from $K_0(C^*(G; \sigma))$ to $\Re$ induced by $\varphi_t$
is locally constant. Denote by $X$ the image of this map. Then
$X$ has cardinality at most the number of components of $\widehat{H_{\sigma}}$, i.e.,
$|(H_{\sigma})_{\tor}|$.
Also, $X$ contains all the extremal points of $T(C^*(G; \sigma))_{K_0}$ and hence
the closed convex hull of $X$ is $T(C^*(G; \sigma))_{K_0}$ by the Krein-Milman theorem \cite[Theorem V.7.4]{Conway}.
Evaluating elements of $X$ at minimal
projections in the center of $C^*(G; \sigma)$ we see that the closed convex hull of $X$
is a simplex of dimension $|H_{\sigma}|-1$ with vertex set $X$.
This proves the first assertion of Lemma~\ref{states:lemma}.

Note that if $_AE_B$ is a strong Morita equivalence bimodule for two unital
$C^*$-algebras $A$ and $B$, and $K(J), J, Y$ are closed two-sided ideals and a submodule
as in Subsection~\ref{CME:subsec}, then the diagram
\begin{equation*}
\xymatrix{
K_0(A) \ar[d] \ar[r] &  K_0(B) \ar[d] \\
K_0(A/K(J)) \ar[r] & K_0(B/J) &
}
\end{equation*}
commutes, where the horizonal isomorphisms are induced by $_AE_B$ and $_{A/K(J)}E/Y_{B/J}$ respectively and
the vertical homomorphisms are induced by the $C^*$-algebra quotient maps.
Using Theorem~\ref{CME to tensor2:thm} one sees that the homomorphism
from $K_0(C^*(G;\sigma))$ to $K_0(C^*(G; \sigma)/I_t)$ induced by the $C^*$-algebra quotient map
is surjective for every maximal ideal $t$ of the center of $C^*(G;\sigma)$.
Thus the image of $K_0(C^*(G;\sigma))$ under any element of $X$ is the image of $K_0(C^*(G/H_{\sigma}; \sigma'))$
under the element of $T(C^*(G/H_{\sigma}; \sigma'))_{K_0}$ induced by $\varphi$.
This finishes the proof of Lemma~\ref{states:lemma}.
\end{proof}

Let $\omega_{\sigma}$ be any vertex of $T(C^*(G; \sigma))_{K_0}$ in Lemma~\ref{states:lemma}.
The following result is a generalization of
Proposition~\ref{ordered K-0:prop}.

\begin{proposition} \label{ordered K-0 2:prop}
Let $\sigma_j$ be a $2$-cocycle on a finitely generated abelian group $G_j$ for $j=1,2$.
Then  $C^*(G_1; \sigma_1)$ and $C^*(G_2; \sigma_2)$
have isomorphic ordered $K_0$-groups if and only if
$\omega_{\sigma_2}(K_0(C^*(G_2; \sigma_2)))=\mu
\omega_{\sigma_1}(K_0(C^*(G_1; \sigma_1)))$ for some real number
$\mu>0$, $|(H_{\sigma_1})_{\tor}|=|(H_{\sigma_2})_{\tor}|$, and
$\rank(G_1)=\rank(G_2)$ or $\rank(G_1)+\rank(G_2)=1$.
\end{proposition}
\begin{proof} Let us prove first the ``only if'' part.
An argument similar to that in the proof of the ``only if'' part
of Proposition~\ref{ordered K-0:prop}
shows that $\omega_{\sigma_2}(K_0(C^*(G_2; \sigma_2)))=\mu
\omega_{\sigma_1}(K_0(C^*(G_1; \sigma_1)))$ for some real number
$\mu>0$. By Lemma~\ref{states:lemma} the set of positive homomorphisms
from $(K_0(C^*(G_j; \sigma_j))_+, K_0(C^*(G_j; \sigma_j)))$ to $(\Re_+, \Re)$ is a cone
of dimension $|(H_{\sigma_j})_{\tor}|$. Thus,
$|(H_{\sigma_1})_{\tor}|=|(H_{\sigma_2})_{\tor}|$.
Comparing the ranks of the $K_0$-groups, by Corollary~\ref{K-groups:coro} we obtain
$\rank(G_1)=\rank(G_2)$ or $\rank(G_1)+\rank(G_2)=1$.

For the ``if'' part, by Theorem~\ref{CME to tensor2:thm} and the ``only if'' part
we may assume that $G_j$ and $\sigma_j$ have the same properties as $G'$ and $\sigma'$ of Theorem~\ref{CME to tensor2:thm}.
Then an argument similar to that in the proof of the ``if'' part
of Proposition~\ref{ordered K-0:prop} completes the proof.
\end{proof}

An argument similar to that in the proof of Lemma~\ref{K-0 for simple:lemma}
establishes the following generalization of
Lemma~\ref{K-0 for simple:lemma}:

\begin{lemma} \label{K-0 for simple 2:lemma}
Let $\sigma_j$ be a skew-symmetric bicharacter on a finitely generated abelian group $G_j$ for $j=1,2$ such that
$G_j=G''_j\oplus G'_j$, $\sigma_j(G''_j, G_j)=1$, and
$G'_j$ is torsion-free.
Denote by $\sigma'_j$ the restriction of $\sigma_j$ to $G'_j$.
Suppose that $\rank(G'_1)=\rank(G'_2)$.
If the
ordered $K_0$-groups of $C^*(G_1;\sigma_1)$ and $C^*(G_2;\sigma_2)$ are
isomorphic,  then so also are those of
$C^*(G'_1; \sigma'_1)$ and $C^*(G'_2;\sigma'_2)$.
\end{lemma}

Now the proof of Theorem~\ref{Morita:thm} in the case of
noncommutative tori in Section~\ref{Morita:sec} extends verbatim to the general case of the
twisted group $C^*$-algebras of
arbitrary finitely generated abelian groups.

\begin{remark} \label{Morita quotient:remark}
For a $2$-cocycle $\sigma$ on a finitely generated abelian group $G$,
one can check easily that the maximal ideals of $C^*(G;\sigma)$
are exactly those $I_t$'s in the proof of Lemma~\ref{states:lemma}, and
hence that the simple quotient algebras of $C^*(G;\sigma)$ are all isomorphic
to $C^*(G/H_{\sigma};\sigma')$ therein. The proof of
Theorem~\ref{CME to tensor2:thm}
actually shows that two twisted group $C^*$-algebras of
finitely generated abelian groups are strongly Morita equivalent if and
only if they have isomorphic centers and their simple quotient
algebras have isomorphic ordered $K_0$-groups.
\end{remark}

\end{document}